\documentclass{amsart}
\usepackage{amssymb}
\title{The complexity of fuzzy logic}
\author{Martin Goldstern}

\newcommand{\<}{\lessdot}
\def\itm#1 {\item[{(#1)}]}
\newcommand{\R}{{\mathcal{R}}}

\newcommand{\LL}{{\mathcal L}}
\newcommand{\M}{{\mathcal M}}

\newcommand\N{{\mathbb N}}
\newcommand{\logand}{\&}
\newcommand{\limpl}{\to}
\newcommand{\liff}{\leftrightarrow}

\newcommand{\on}{{|}}
\newcommand{\wert}[1]{\left\| #1 \right\|}
\newcommand{\class}{{\mathcal L}^{\text{class}}}
\newcommand{\pc}{\class_\sigma}
\newcommand{\Add}{{\it Add}}
\newcommand{\Mul}{{\it Mul}}

\def\xx#1 {%
\newtheorem{#1}[thm]{#1}}
\def\yy#1 {%
\newenvironment{#1}{\begin{roman-#1}\rm}{\end{roman-#1}}
\newtheorem{roman-#1}[thm]{#1}}

\yy Fact
\xx Theorem
\xx Lemma
\yy Definition
\yy Remark
\yy Notation
\yy {Fact and Definition}
\yy Acknowledgements
\xx {Main Claim}
\yy Claim
\xx Observation
\xx Explanation

\begin{document}
\nocite{Hajek2}

\begin{abstract}
{We show  that the set of 
valid formulas in {\L}ukasiewicz predicate logic is a complete 
$\Pi^0_2$ set.   We also show that the classically valid formulas
are exactly those formulas in the classical language whose fuzzy value
is $1/2$. 
}
\end{abstract}
\date{July 1997}

\maketitle

{\L}ukasiewicz' infinite valued logic can be seen as a 
particular ``implementation'' of fuzzy logic.    The set of 
possible ``truth values'' (or, in another interpretation, 
degrees of certainty) is the real interval $[0,1]$.   Minimum,
maximum, and truncated addition are the basic operations. 

It is well known that the {\em propositional} fragment  
version of 
{\L}ukasiewicz logic is decidable.

The exact complexity of {\L}ukasiewicz {\em predicate} logic was a
more diffcult problem.  For an upper bound, it is known that the set of
valid formulas in this logic is a $\Pi^0_2$.

\nocite{Hay}
(An explicit $\Pi^0_2$ representation can be found through the 
axiomatisation of Novak and  Pavelka.   See \cite{Hajek2} for
references.)

For a lower bound, Scarpellini \cite{Sc} showed that the set of 
valid formulas is
not r.e., and in fact 
 $\Pi^0_1$-hard.  He also remarks in a footnote that this set is not
$\Sigma^0_2$, either.    

In his unpublished thesis \cite{Ragaz}, Mathias Ragaz showed that the set of
valid formulas in {\L}ukasiewicz predicate logic is actually
$\Pi^0_2$-complete.  The proof of this theorem that we give here was 
found independently.

Furthermore, we show that 
if we restrict our attention to classical
formulas, the classically valid formulas are exactly those formulas
which have value $\ge 1/2 $ in every fuzzy model.

\section{Definitions}

For the reader's convenience, we recall the syntax and
semantics of {\L}ukasiewicz' logic.

\subsection{The natural MV-algebra on $[0,1]$:}

We define a unary operation $\lnot$ and binary operations 
$\limpl$, $\wedge$, $\vee$, $\&$, $+$ on the unit interval 
$[0,1]$ as follows:
$\lnot r = 1-r$, 
 $r \vee s = \max(r,s)$, 
$r+s = \min(1, r+s)$ (where $+$ on the right side is ``true''
addition), and $\wedge $ and $\&$ are dual to $\vee$ and $+$:
$r \wedge s = \min(r,s)$, $1 - (r \& s)  = (1 - r) + (1 - s)$. 
We let $r \limpl s = (\lnot r ) + s $.

($\&$ and $+$ are called the ``strict conjunction'' and 
``strict disjunction''.)



\subsection{\bf Propositional \L-formulas and assignments:}

Propositional {\L}ukasiewicz logic uses the connectives $\limpl$,
$\lnot$, $\wedge$, $\vee$, $\&$ and $+$ (sometimes also written as
$\underline{\vee}$).    Formulas are built in the usual way from 
an infinite set of propositional variables. 

An assignment is a map  from the variables into $[0,1]$. 
 Each assignment  $s$ naturally induces  a map $\bar s$
 from the set of all formulas into $[0,1]$:  $\bar s$ extends $s$, 
$\bar s(\varphi  \logand \psi ) = \bar s(\varphi)
 \logand \bar s(\psi)$, etc. 

We let $\wert{\varphi } = \inf \{s(\varphi): 
\hbox{$s$ an assignment}\}$, and we 
call a formula $\varphi$ an \L-tautology iff $\bar
s(\varphi)=1$ for all assignments $s$, i.e., if $\wert{\varphi}=1$.


\bigskip

 \subsection{Predicate \L-logic and models}

In \L ukasiewicz predicate logic, we have in addition to the
connectives also quantifiers $\forall$ and $\exists$.

Let $\sigma $ be a (finite) set of relation symbols and constant
symbols, with an ``arity'' attached  to each relation symbol.  

The set of terms and formulas in the language $\LL_\sigma$ for fuzzy
predicate calculus over the signature $\sigma$ is defined in the usual
way, starting from constants and variables. 


A `model' $\M$ for the language $\LL_\sigma$ is given by 
a nonempty set $M$, together with interpretations $c^\M$ and $R^\M$ of
all constants and relation symbols.  Each $c^\M$ is an element of
$\M$, and if $R$ is a $k$-ary relation symbol, then $R^\M$ is a map
from $M^k$ into $[0,1]$. 

$\LL_\sigma(M)$ is the language $\LL_\sigma$ expanded by special constant
symbols $c_m$ for every $m\in M$. We will not notationally distinguish 
between $m $ and $c_m$, and we require that always $c_m^\M=m$. 

An $M$-formula is a formula in $\LL(M)$. 

 $\wert \varphi^\M \in [0,1]$ is defined by induction for all closed
$\M$-formulas in the natural way:  $\wert{R(c)} = R^\M(c^\M)$ for any
constant symbol $c$, $\wert{\lnot \varphi} = \lnot \wert{\varphi}$, 
$\wert{ \forall x \, \varphi(x) } = \inf\limits_{a\in M}
\wert{\varphi(a)}^\M$, etc. 
\\
We let $\wert{\varphi}$ be the infimum over $\wert{\varphi}^\M$, 
taken over all fuzzy models $\M$.

\subsection{\L-validities and \L-tautologies}

We say that a closed formula  $\varphi$ is \L-valid if 
$\wert{\varphi} = 1$, i.e., if $\wert{\varphi}^\M = 1 $ for all 
models $\M$.

We call a formula $\varphi$ in $\LL_\sigma$ (or even in $\LL_\sigma(M)$)
an  {\em \L-tautology} iff there is 
an \L-tautology $\chi$ in the propositional fragment of 
\L ukasiewicz logic and a homomorphism $h$ that assigns formulas 
in $\LL_\sigma$ to propositional variables such that $\varphi = h(\chi)$. 

Thus, the set of \L-tautologies is a (decidable, proper) subset of the
set of all        \L-validities. 


\section{The complexity of \L-validities}

We now proceed to prove the following theorem.

\begin{Theorem}\label{thm}
Let $\sigma$ be a sufficiently rich  (relational) signature, 
$\LL_\sigma$ the set of formulas for this signature in 
{\L}ukasiewicz logic.   
Then the set of formulas in $\LL_\sigma$ which is ``valid'' is 
 is $\Pi^0_2$-complete. 
\end{Theorem}

\begin{Definition} \label{04-21}
Let $\R \subseteq \N \times  \N $ be a primitive recursive
relation such that 
 $$A:=\{m: \exists^\infty n \, (m,n)\in \R\}$$ is
$\Pi^0_2$-complete. 

If $m\notin A$, let $f(m)$ be a positive  natural number 
such that $\forall n\ge f(m): (m,n) \notin \R$. 

\end{Definition}

\begin{Definition}

We start with an extended (purely relational) language $\LL_{\sigma_0}$ of
Peano arithmetic. The signature $\sigma_0$ contains  
constant symbols  $0$, $1$, binary relation symbols  $<$ and 
$\<$ (read:  ``is direct successor of'') 
 and finitely many relation symbols (such as  ternary 
$\Add$ and $\Mul$) intended to code various
primitive recursive relations, among them a symbol $R$ intended to
be interpreted by  $\R$.

Our full signature $\sigma $ will have in addition to the 
symbols in $\sigma_0$ two
extra predicates 
(which are intended to be fuzzy): 
A binary relation  $Q$ and unary relation $P$.  
\end{Definition}

\begin{Definition}
Let $\varphi_0$ be an axiom similar to Robinson's Q, coding enough
primitive recursive definitions such that in classical logic the
following are derivable: 
\begin{enumerate}
\item[$(+)$]  Whenever $(m,n)\in \R$, then $$\varphi_0 \vdash 
0 \< 1 \< x_2\<\cdots \< x_m\  \wedge \ 
0 \< 1 \< y_2 \< \cdots \< y_n \limpl R(x_m, y_n)$$
\item[$(-)$]  Whenever $(m,n)\notin \R$, then $$\varphi_0 \vdash 
0 \< 1 \< x_2\<\cdots \< x_m\  \wedge \ 
0 \< 1 \< y_2\<\cdots \< y_n \limpl \lnot R(x_m, y_n)$$
Moreover, we assume $\varphi_0 \vdash  \forall x \, \exists y \,\, x \< y$
\end{enumerate}
\end{Definition}

\begin{Definition}
The formula $\varphi_1 $ is the universal closure of 
$$ R(x,y) \logand Q(y,y)  \limpl P(x)$$
The formula $\varphi_2 $ is the universal closure of 
$$  Q(1,y) \liff  \lnot Q(y,y)$$
The formula $\varphi_3 $ is the universal closure of 
$$  x \< x' \ \limpl \ [( Q(x,y){+} Q(1,y) ) \liff Q(x',y)] $$
\end{Definition}

\begin{Definition}
The {\em standard model} $\N$  is defined as follows:
\begin{itemize}
\item[$*$] The universe of $\N$ is the set of natural numbers.
\item[$*$] $\Add$, $\Mul$, $<$, $\<$, \dots, $R$ are interpreted naturally.
\item[$*$] $\wert{Q(m,n)}^\N = \min(1, m/(n+1))$ for all $m,n\in \N$.
\item[$*$] $\wert{P(m)}^\N = 1 $ if $m\in A$
\item[$*$] $\wert{P(m)}^\N = 1 - \frac{1}{f(m)}$ if $m\notin A$
(where $f$ is the function defined in \ref{04-21}). 
\end{itemize}
\end{Definition}

\begin{Fact}
$\wert{\varphi_1 \wedge \varphi_2 \wedge \varphi_3 \wedge \varphi_4}^\N=1$.
\end{Fact}

\begin{Definition}

For each formula $\varphi$, let $\varepsilon_\varphi$ be the universal
closure of $\varphi \wedge  \lnot \varphi$. 

Let $\varepsilon$ 
be the disjunction over all $\varepsilon_\varphi$,
where $\varphi$ runs over all subformulas of $\varphi_0$, $\varphi_1$,
$\varphi_2$, $\varphi_3$ {\bf except} the formulas $Q(x,y)$ and $P(x)$. 
\end{Definition}

  So $\varepsilon$  measures 
how close a structure (disregarding $Q$ and $P$) comes
to being a crisp model.  In the structures that we are interested in 
the value of $\varepsilon $ will be $0$ (or at least close to 0).

Note that  for any formula $\varphi(x,y,\ldots)$ that was used in
the definition of $\varepsilon$, any structure  $\M$,
and any $a,b\ldots$ in $\M$ we have: If $ e := \wert{\varepsilon}^\M$, then 
$\wert{\varphi(a,b,\ldots)}^\M \in [0,e] \cup [1-e, 1]$.

\begin{Definition}
For each $m>3$, let $\psi_m$ be the universal closure of the 
following formula: 
$$ \varphi_0 \logand \varphi_1 \logand \varphi_2 \logand \varphi _ 3
\logand (0 \< 1\<x_2 \< \cdots \< x_m )\ 
\limpl
\ 
P(x_m) +  10\varepsilon $$
where $2 \varepsilon = \varepsilon + \varepsilon$, $3 \varepsilon  =
2 \varepsilon  + \varepsilon $, etc.   The formula $0\< 1 \< \cdots$
is any conjunction (it does not matter if sharp or not) of the formulas
$0\<1$, $1\< x_2$, \dots 
\end{Definition}

\begin{Explanation}
In a fuzzy model where $\varphi_0$ is  true (or at least ``sufficiently
true'', i.e., has a truth value close to 1), the formula
$\exists x_2 \cdots \exists x_m\,\, 0 \< 1 \< x_2 \< \cdots \< x_m $
says that $x_m$ is similar to the number $m$. 

In a model where $\varphi_2
\logand \varphi_3$ is (sufficiently) true, the formula $Q(y,y)$
expresses the fact that $y$ is ``infinite'', or at least ``large''. 

In a model of $\varphi_1$, $P(x)$ says that $R(x,y)$ holds for some
large $y$. 

Hence, $\psi_m$ says  that an object that has properties similar to
the number $m$ is in a set that is similar to $A$. 

This is only an approximation, of course.  We wil lnow show that this
approximation is good enough for our purposes. 
\end{Explanation}

\begin{Main Claim}
For all $m>3$: 
 $m \in A $ iff $\wert{ \psi_m} = 1 $. 
\end{Main Claim}

Clearly this claim will imply our theorem. 

\begin{proof}[Proof of the main claim, part 1]
First assume that $m $ is not in $A$.  So for all $n\ge f(m)$,
$(m,n)\notin \R$.  
Evaluate $\psi_m$ in the standard model $\N$: 

Clearly $\varepsilon$ will have value $0$. 

Instantiate $\psi_m$ with $x_i=i$. 
 Now the antecedent  has value 1, whereas the
consequent has value $<1$, so $\wert{\psi_m} < 1$, and we are done. 
\end{proof}
\bigskip

\begin{proof}[Proof of the main claim, part 2]
Now assume that $m\in A$. Fix a small number ${\delta}$ (in particular
${\delta} < 1/m$).  We will show
that in all fuzzy structures,  all instances of $\psi_m$ have a value $\ge
1-{\delta}$. 

So fix a structure $\M$, and let $a_2$, \dots, $a_m$ be elements of $\M$. 
Write $e$ for the value $\wert{\varepsilon}^\M$. 

   Let $\psi'_m$ be the instance of $\psi_m$ obtained by
substituting $a_i$ for $x_i$, i.e., 
$$ \varphi_m' \ = \ 
\varphi_0 \logand \varphi_1 \logand \varphi_2 \logand \varphi _ 3
\logand (0 \< 1\<a_2 \< \cdots \< a_m )\ 
\limpl
\ 
P(a_m) +  10\varepsilon $$
If  ${e} \ge 0.1$, then clearly the value of $\psi '_m$ is
1.   So assume that ${e} < 0.1$.    Hence all fuzzy
relations associated with symbols in the language $\LL_0$ will always
take values $<0.1$ or $>0.9$. 

Let $\M'$ be the structure obtained from $\M$ by keeping the values of
$Q$ and $P$, but rounding all other basic relations to 0 or 1. 

We will write $r \sim s$ if $|r-s|<0.1$. 

\begin{Claim}
Whenever $\chi$ is a closed $M$-formula which is a 
substitution instance of some subformula of
 $\psi'$ (except for $Q( \cdot, \cdot)$ 
and $P(\cdot)$), we have $\wert{\chi}^\M \sim \wert{\chi}^{\M'}$.
\end{Claim}
The proof of the claim is by induction on the complexity of
$\chi$, using the fact that [by the definition of $\varepsilon$]
we must have $\wert{\chi}^\M \in [0,e]\cup [1-e,1]$, and $e+e<1-e$.

\bigskip

We continue the proof of part 2 of the main claim:

If any of the formulas $\varphi_0$, $\varphi_1$, $\varphi_2$,
$\varphi_3$, $0 \< 1$, $1\< a_2$, $a_i \< a_{i+1} $
 has (in $\M$) a value $\le {e}$, then 
$\wert{\psi'_m}^\M = 1$.  So assume that all these values are $\ge 1-
{e}$.   By the above claim, $\wert{\varphi_0}^{\M'} = 1$. 

 Choose some large enough $n$ (specifically, let  $n >
1/{\delta}$, and if ${e} >0$ then also $n> 1/{e}$),
such that $(m,n)\in \R$. 

Working in $\M'$, we can use the fact that $\M'\on \LL_0$ satisfies 
$\varphi_0$ in the classical sense:
\begin{enumerate}
\item We can find $a_{m+1}$, \dots, $a_n$ such that
for all $i< n$: $\M' \models a_i \< a_{i+1}$. 
\item  $\M' \models R(a_m,a_n)$. 
\item  Therefore,  $\wert{R(a_m, a_n)}^\M \ge 1 - {e}$. 
\end{enumerate}

We now claim that 
$$ (**) \qquad \wert{Q(a_n, a_n)}^\M \ge \min(1-4e, 1-\delta)$$
Using $(**) $ and the fact that 
$\wert{ R(a_m,a_n) \logand Q(a_n,a_n)  \limpl P(a_m)} \ge 
\wert{\varphi_1}^\M \ge 1 - {e} $, we can then conclude 
$$\wert{P(a_m)}^\M \ge \min(1 - 6 {e}, 1-\delta-2e), $$
hence $\wert{\psi'_m }^\M \ge 1-\delta$, so $(**)$ would 
finish the proof of the main claim. 

 \medskip

So it remains to prove $(**)$. 

Let us abbreviate $q_1:=Q^\M(1,a_n)$ and 
$q_k:=Q^\M(a_k, a_n)$ for $k=2,\ldots, n$.

Then we have for $k=1,\ldots, n$: 
$$ \wert{ a_k \< a_{k+1} \ \limpl \ [( Q(a_k,a_n){+} Q(1,a_n) ) \limpl
Q(a_{k+1},a_n)]}^\M \ge 1-{e}. $$
Since also $\wert{a_k \< a_{k+1}}^\M \ge 1-{e}$, we get 
$$ \wert{ ( Q(a_k,a_n){+} Q(1,a_n) ) \limpl
Q(a_{k+1},a_n)}^\M \ge 1-{2e},  $$
so 
$q_{k+1} \ge (q_k+q_1) -2e$.     Using induction on $k$  we can show 
for all $k$: 
\begin{itemize}
\item[$(*1)_k$] $ q_k \ge \min(k\cdot (q_1-2e), 1-2e)$.
\end{itemize}
In particular, we get 
\begin{itemize}
\item[$(*1)$] $ q_n \ge \min(n\cdot (q_1-2e), 1-2e)$.
\end{itemize}
We also have $\wert{\lnot Q(1,a_n) \limpl Q(a_n,a_n)}^\M\ge 1-e$, 
so 
\begin{itemize}
\item[$(*2)$]
$q_n \ge 1-q_1 -e$.
\end{itemize}

\smallskip

In the proof   of $(**)$ we distinguish 3 cases:
\begin{enumerate}
\item $e=0$.
\item $e>0$, $Q^\M(1,a_m) \le 3e$.
\item $e>0$, $Q^\M(1,a_m) > 3e$.
\end{enumerate}

\noindent{\bf Case 1:}  ${e} = 0$. \\
Assume $q_n < 1$.  By $(*2)$, $q_1 \ge 1-q_n$.  By $(*1)$, 
$q_n \ge n\cdot q_1 \ge n\cdot (1-q_n)$.  So 
$q_n \ge n/(n+1) > 1- \delta$. 

\noindent{\bf Case 2:}  ${e} >0$,  $q_1 \le 3e$.\\   
By $(*2)$, we have $q_n \ge  1-4e$. 

\noindent{\bf Case 3:}  ${e} >0$,  $q_1> 3e$.\\   
Recall that $n$ was chosen so large that $n\cdot e > 1$, so 
$n \cdot (q_1 - 2e) > 1$, so 
$(*1)$ implies $q_n \ge 1-2e$.

\end{proof}


\section{A restricted language}

In the previous section we have made good use of the fact that 
addition is hardwired into the semantics of our particular brand 
of fuzzy logic.   In this section we show that if we restrict the 
language to the lattice operations $\vee$, $\wedge$, together with 
$\lnot$, then the computation of $\wert{\varphi}$ can be 
reduced to the problem of deciding classically valid formulas, and
 conversely. 

Let $\class$ be the ``classical'' propositional 
language, using only the connectives
$\wedge$, $\vee$ and $\lnot$, and let $\pc$ be the  
classical predicate language over the  signature $\sigma$. 

We mention some easy (and well-known) fact:

\begin{Observation}\label{04-23}
If $\varphi$ is a propositional formula in $\class$, then
\begin{enumerate}
\itm a For any $\varphi \in \class$,
 $\wert{\varphi} \le \frac{1}{2}$. 
\itm b If $\varphi$ is a classical tautology, then $s(\varphi)\ge
0.5$ under any assignment $s$. 
\itm c If $\varphi$ is not a classical tautology, then  $\wert{\varphi} = 0$. 
\itm d  The following are equivalent for $\varphi \in \class$: 
\begin{itemize}
\item $\varphi$ is a classical tautology
\item $\wert{\varphi} = \frac{1}{2}$
\item $\wert{\varphi} > 0 $. 
\item $\wert{\lnot \varphi \limpl \varphi} = 1 $.
\item $\wert{p \wedge \lnot p \limpl \varphi}$, where $p$ is 
	any propositional variable not appearing in $\varphi$. 
\end{itemize}
\itm e The set of fuzzy tautologies is co-NP-complete. 
\end{enumerate}
\end{Observation}
\begin{proof}

(a): Assign the value $\frac 1 2$ to every propositional variable. 

(c) is clear, and (d) is a reformulation of (a)--(c). (e) follows 
from (d).

It remains to prove  (b): 
\\
Let us call a formula $\varphi$ a ``literal'' if $\varphi = p$ or 
$\varphi = \lnot p$ for some   propositional variable $p$. 
A ``clause'' will be a nonempty conjunction of literals.   
We say that $\psi$ is in
``normal form'' if $\psi$ is a nonempty  disjunction of clauses.

Using the distributive law, as well as de Morgan's laws and cancelling
of double negations, we can find a formula $\psi$ which is equivalent
(classically as well as in propositional \L-logic) 
to $\varphi$ and is in normal form. 
\\
Now, if $\psi   $ is a classical tautology, 
then  each clause of $\psi$ contains some variable
$p$ both in negated and unnegated form. Hence, under any fuzzy
assignment $s$, $s(\psi ) \ge 0.5$.
\end{proof}

\begin{Theorem}
Let $\varphi$ be a formula in Lukasiewicz predicate logic 
which in  $\pc$.  Then:
\begin{enumerate}
\itm a There is a model $\M$ such that $\wert\varphi^\M = 0.5$. 
\itm b If $\varphi$ is classically valid, then 
	$\wert \varphi = 0.5$.
\itm c If $\varphi$ is not classically valid, then 
	$\wert \varphi = 0$. 
\end{enumerate}
\end{Theorem}

\begin{proof}
(a) and (c) are clear. 

Instead of (b), it is enough to show the following:
\begin{quote}
If $\lnot \varphi$ is
classically valid, then 
$\wert{\varphi}\le 0.5$. 
\end{quote}

Without loss of generality we may assume that $\varphi$ is in 
prenex form (since a transformation to prenex form preserves the
classical truth value as well as $\wert{\varphi}$), say 
$\varphi = (\forall x_1)(\exists y_1) (\forall x_2) \cdots
(\exists y_n) \psi(x_1, y_1, \ldots, x_n, y_n)$, where $\psi$ 
is quantifier-free.  

Let $\bar \psi$ be the Skolemization of $\psi$. That is, let 
$g_i$ be an $i$-ary function symbol for $i=1,\ldots, n$, and let 
$$ \bar \psi(x_1, \ldots, x_n)\  := \ 
  \psi(x_1, g_1(x_1), x_2, g_2(x_1,
x_2), \ldots, x_n, g_n(x_1, \ldots , x_n))$$
and let $\bar \varphi = \forall x_1 \, \forall x_2 \, \cdots \forall
x_n \,\, \bar \psi(x_1, \ldots, x_n)$.

Since $\varphi$ is a classical contradiction, also 
$\bar \varphi$ has
no classical model.   By Herbrand's theorem, there is a finite
conjunction 
$\bar \psi_1 \wedge \cdots \wedge \bar \psi_k$
of closed instances of $ \bar \psi$ which is  a  
contradiction in the sense of classical propositional logic. 

[More formally, there is a propositional formula $\chi \in \class$
which is a classical propositional contradiction,
 and a homomorphism $h$ such that 
$\bar \psi_1 \wedge \cdots \wedge \bar \psi_k = h(\chi)$.]

Let $\bar \psi_i = \bar \psi(t^i_1, \ldots, t^i_n)$, where all the 
$t^i_j$ are closed terms (involving the function symbols $g_1$, \dots,
$g_n$ and some constant symbols from the original signature.

Now let $\M$ be a fuzzy model, and assume, towards a contradiction,
than $\wert{\varphi}^\M > 0.5$.   For each $m_1 \in M$ choose 
$f_1(m_1)\in M$ such that 
$$\wert{(\forall  x_2) (\exists y_2)\cdots (\forall x_n)(\exists y_n) 
\psi(m_1, f_1(m_1), x_2, y_2, \ldots, x_n, y_n)}^\M > \frac{1}{2}$$

Now 
for each $m_1, m_2$ in $M$ choose
$f_2(m_1, m_2)\in M$ such that 
$$\wert{(\forall x_3) \cdots (\exists y_n) \psi(m_1, f_1(m_1), m_2, 
f_2(m_1, m_2), \ldots, y_n)}^\M > \frac{1}{2} $$
and continue by induction. 
We thus get functions $f_1, \ldots, f_n$ such that for any 
$a_1,\ldots, a_n \in M$: 
$$\wert{\psi(a_1,f_1(a_1), a_2, f_2(a_1,a_2), \ldots,
a_n, f_n(a_1,\ldots, a_n))}^\M > \frac{1}{2} $$

Now $\bar\M := (M, f_1, \ldots, f_n, c^\M: c \in \sigma)$ 
is a (classical) structure for the 
signature $(g_1, \ldots, g_n, c:{c\in \sigma})$.   

For any (quantifier-free) closed $\bar \M$-formula
 $\chi$ let $\chi^{\bar \M}$ be
obtained by replacing each atomic subformula $P(t_1, \ldots)$ 
by $P(t_1^{\bar \M}, \ldots)$, where $t^{\bar \M}$ is the value of the
closed  term $t$. 

Now recall that $\bar \psi_1 \logand \cdots $ was a classical
contradiction.  So also $(\bar \psi_1 \logand \cdots )^{\bar \M}$ is a
classical contradiction.  Note that this formula does not contain any
function symbols any more, so we can compute its value in our fuzzy
structure  $\M$. By \ref{04-23}, this value is at most $1/2$, so wlog
$\wert{ \bar \psi_1^{\bar \M}}^\M \le 1/2$. 

Let $a_j := {t^1_j}^{\bar \M}$. 
Now
$$ 
\begin{array}{rcl}
\wert{ \bar \psi_1^{\bar \M}}^\M &=&
\wert{ \bar \psi(t^1_1, \ldots, t^1_n)^{\bar \M}}^\M = \\
&=& \wert{ \psi(t^1_1, g_1(t^1_1), t^1_2, \ldots, t^1_n, 
g_n(t^1_1, \ldots, t^1_n)) ^{\bar \M}}^\M = \\
&=&\wert{ \psi(a_1, f_1(a_1), a_2,  \ldots) }^\M > 1/2
\end{array}
$$
a contradiction.

\end{proof}

\vfill
\medskip\hrule\medskip
Please send comments to  {\tt Martin.Goldstern@tuwien.ac.at}. 
\medskip\hrule\medskip
\vfill
\eject
\end{document}